\documentclass[oneside,pdflatex,sn-mathphys-num,numbered]{sn-jnl}
\usepackage[T1]{fontenc}
\usepackage[utf8]{inputenc}
\usepackage{color}
\usepackage{float}
\usepackage{booktabs}
\usepackage{amsmath}
\usepackage{amsthm}
\usepackage{amssymb}

\makeatletter

\floatstyle{ruled}
\newfloat{algorithm}{tbp}{loa}
\providecommand{\algorithmname}{Algorithm}
\floatname{algorithm}{\protect\algorithmname}

\numberwithin{equation}{section}

\@ifundefined{date}{}{\date{}}
\let\origprintabstract\printabstract
\renewcommand{\printabstract}{\centerline{\normalsize November 14, 2026}\vspace{1em}\origprintabstract}

\makeatother

\providecommand\theoremname{Theorem}
\theoremstyle{plain}
\newtheorem{thm}{\protect\theoremname}
\providecommand\definitionname{Definition}
\theoremstyle{definition}
\newtheorem{defn}[thm]{\protect\definitionname}
\providecommand\propositionname{Proposition}
\theoremstyle{plain}
\newtheorem{prop}[thm]{\protect\propositionname}
\providecommand\lemmaname{Lemma}
\newtheorem{lem}[thm]{\protect\lemmaname}
\providecommand\examplename{Example}
\theoremstyle{definition}
\newtheorem{example}[thm]{\protect\examplename}
\providecommand\remarkname{Remark}
\theoremstyle{remark}
\newtheorem{rem}[thm]{\protect\remarkname}
\providecommand\corollaryname{Corollary}
\theoremstyle{plain}
\newtheorem{cor}[thm]{\protect\corollaryname}
\usepackage[english]{babel}
\begin{document}
\title{\textbf{Remotest Operator Control String-Averaging Methods in Feasibility-Seeking}}
\author{Yair Censor\textsuperscript{1}}
\email{yair@math.haifa.ac.il}
\author*{Ariel Nisenbaum\textsuperscript{2*}}
\email{nariel@campus.technion.ac.il}
\author{Simeon Reich\textsuperscript{2}}
\email{sreich@technion.ac.il }
\affil{\textsuperscript{1}Department of Mathematics, University of Haifa,
Mt. Carmel, Haifa 3498838, Israel}
\affil{\textsuperscript{2*}Department of Mathematics, The Technion ---
Israel Institute of Technology, Haifa 3200003, Israel }
\abstract{The string-averaging (SA) method is a versatile algorithmic scheme
for solving feasibility-seeking problems. It partitions the index
set of all constraints into subsets of indices, called strings, allowing
the user to apply operators successively along each string separately,
and thus in parallel. The end-points of all such string operations
are convexly combined to create the next iterate. Due to its flexible
framework, SA is not restricted solely to this procedure; however,
this characteristic applies only to the simultaneous step, while the
sequential one is restricted to operator compositions. To permit alternative
strategies in the latter, a novel algorithmic structure termed the
Skeletal SA is introduced. Although its convergence is not established,
it serves as a foundational \textquotedbl parent\textquotedbl{} structure
that unveils unexamined SA algorithmic structures, as well as other
known methods, obscured by standard compositions. To demonstrate the
benefit of adopting the Skeletal SA framework over the standard SA
scheme, we investigate a new approach, where instead of passing sequentially
over each string via compositions, we employ a remotest-set-control-like
strategy on every string by identifying the index of the operator
farthest from the current iterate, and invoking it. This method favorably
handles the common fixed point problem (CFPP), which requires finding
a point within the intersection of finitely many constraint sets whose
intersection is nonempty. Many iterative methods, including SA ones,
have been proposed to solve the CFPP due to its importance in many
scientific fields. However, remotest operator control SA (ROC-SA)
methods have remained unexamined in this context.}

\maketitle
Keywords: String-averaging; Skeletal string-averaging; Remotest operator
control; Common fixed point problem; Projection methods; Convex feasibility
problem.

\section{Introduction}

A sequential algorithm is a method in which in every iteration all
steps are done one after the other. This is in contrast to a simultaneous
(or parallel) algorithm, that carries out several similar computations
(steps) at once. The string-averaging (SA) algorithmic structure,
first introduced in \cite{Pink-book}, combines these two common schemes.
It handles a family of operators, which needs to be applied sequentially
along the indices of ordered index sets, known as ``strings'', and
the obtained end-points are then combined in some agreed manner to
yield the next iterate of the procedure.

Since its presentation, the notion of SA has been widely studied;
see, e.g., \cite{CensorZaslavski1,CensorZaslavski2,CQmethods,ReichZalas,SegalCensor,kong2019,BMR},
to name but a few. Nevertheless, to the best of our knowledge, prior
work has focused on using SA to expand the class of methods for solving
feasibility-seeking problems or to exploit its ability to handle large-scale
real-world applications in terms of computing time, and have not addressed
its adaptable algorithmic structure. In \cite{GordonANDGordon} a
procedure which possesses similar nature to the SA scheme is proposed.
The sequential step described there carries out a finite number of
Kaczmarz's algorithm iterative steps with the option of using a relaxation
parameter, see \cite[Algorithm 1]{GordonANDGordon}. That is, it involves
an additional action apart from composition. This raises the question
whether the sequential part of the SA algorithmic structure can be
modified to a different ``sequential type'' action, as well. 

Motivated by this question, the Skeletal SA structure is introduced
herein. To demonstrate its potential, we define, investigate and study
``remotest operator control SA (ROC-SA)'' methods, where at every
iteration the operators associated with the most violated constraint
in each string are selected, and then combined through a convex combination
to yield the next iterate. Such a scheme is intended to solve the
common fixed point problem (CFPP); see, e.g., \cite[Section 1.3.3]{Cegielski},
in which for a given finite family of operators $\left(T_{i}\right)_{i\in I}$
in a real Hilbert space $\mathcal{H}$, one seeks to find a point
in the nonempty intersection $F:=\cap_{i\in I}\text{Fix}T_{i}\neq\emptyset$,
where $\text{Fix}T_{i}$ is the fixed point set of $T_{i}$ for every
$i\in I$; see, e.g., \cite[Equation 1.68]{BauschkeCombettes}.

This well-established problem has widespread use because numerous
problems in various branches of mathematical and physical sciences
can be recast into it, thus, increasing the importance of this problem,
as well as studying algorithms for solving it; see, e.g., \cite{BauschkeCombettes_paper,Combettes,Halpern,HLWB,SegalCensor,Crombez1}
and the excellent references in \cite{BauschkeBorwein}.

Besides the fact that the remotest set control framework plays an
important role in fixed point theory and convex feasibility problems
in a Hilbert space, see, for instance, \cite{ZaslavskiGreedy,GubinPolyakRaik,Agmon,MotzkinSchoenberg}
and \cite[Section 5.6]{Cegielski}, our choice of studying an operator
based variant of it stems from the growing interest in recent years
in the so-called ``greedy'' search algorithms. Although some of
these algorithms inherently feature a control strategy, this specific
terminology is rarely used in their context; see, for example, \cite{GreedyAlgorithmVIP}.
However, since the choice of the ``farthest set (or operator)''
is a greedy one, these iterative procedures can indeed be characterized
as greedy methods.

This paper is organized as follows: After introducing basic terms
and essential tools in Section \ref{sec:Preliminaries}, the Skeletal
SA framework and several of its related variations are presented in
Section \ref{sec:The-string-averaging-algorithmic}. Section \ref{sec:Main-Results}
is devoted to establishing the main convergence results concerning
our ROC-SA scheme. Some of our methods there are labeled as “static”
because they employ a fixed, pre-determined set of strings, while
others are labeled as “dynamic” because they allow the selection of
strings to vary dynamically across iterations. Concluding remarks
and future research directions are outlined in Section \ref{sec:Concluding-Comments}.

\section{Preliminaries}\label{sec:Preliminaries}

We denote by $\mathcal{H}$ a real Hilbert space with the inner product
$\left\langle \cdot\,,\,\cdot\right\rangle $ and induced norm $\Vert\cdot\Vert$.
Denote by $\mathbb{N}$ the set of all natural numbers including zero.
Let $D\subseteq\mathcal{H}$. If $D$ is a closed and convex set,
then for every $x\in\mathcal{H}$ there exists a unique point $y\in D$
such that 
\begin{equation}
\Vert x-y\Vert=d\left(x,D\right):=\text{inf}_{z\in D}\Vert x-z\Vert.\label{eq: inf-distance}
\end{equation}
Such a point $y$ is called ``the projection of $x$ onto $D$''
and denoted by $P_{D}\left(x\right)$ . For a finite family of closed
and convex sets $\left(D_{i}\right)_{i\in I}$, we denote the projection
onto each $D_{i}$ by $P_{i}:=P_{D_{i}}$.

It is possible to characterized $P_{D}$ by 
\begin{equation}
\left\langle x-P_{D}\left(x\right),z-P_{D}\left(x\right)\right\rangle \leq0,\label{eq: characterization of projection}
\end{equation}
for every $z\in D$; see, e.g., \cite[Proposition 3.5]{ReichBook}.
(\ref{eq: characterization of projection}) combined with \cite[Lemma 1.2.5]{Cegielski}
yield a simple, yet useful, tool, which states that, for every $z\in D$,
\begin{equation}
\left\Vert z-P_{D}\left(x\right)\right\Vert ^{2}\leq\left\Vert z-x\right\Vert ^{2}-\left\Vert P_{D}\left(x\right)-x\right\Vert ^{2}.\label{eq: characterization and Lemma 1.2.5}
\end{equation}

Another property of the projection operator is its bounded regularity;
see, e.g., \cite[Example 3.5]{CegielskiZalas}. Boundedly regular
operators were first proposed by Cegielski and Zalas in \cite[Definition 16]{CegielskiZalas1}
under the name ``approximately shrinking'', due to their relation
to quasi-shrinking operators \cite[Section 3]{quasiShrinking}. They
belong to the class of quasi-nonexpansive operators, whose definition
is presented next, followed by that of the former.

The first definition is written in a similar manner to \cite[Definition 2.15]{ReichZalas},
adapted to our notations and needs. Before defining those, we note
that for every operator $T:D\rightarrow D$, we denote by $\text{Fix}T:=\left\{ x\in D\mid T\left(x\right)=x\right\} $
the fixed point set of $T$.
\begin{defn}
\label{def:QNE}An operator $T:D\rightarrow\mathcal{H}$ is quasi-nonexpansive
(QNE) if 
\begin{equation}
\left\Vert T\left(x\right)-y\right\Vert \leq\left\Vert x-y\right\Vert ,
\end{equation}
for all $x\in D$ and $y\in\text{Fix}T$.
\end{defn}

\begin{defn}
\label{def:BR operators}A QNE operator $T:\mathcal{H}\rightarrow\mathcal{H}$
is boundedly regular (BR) if for any bounded sequence $\left(x^{k}\right)_{k\in\mathbb{N}}\subseteq\mathcal{H}$,
we have
\begin{equation}
\underset{k\rightarrow\infty}{\text{lim}}\left\Vert T\left(x^{k}\right)-x^{k}\right\Vert =0\:\Longrightarrow\:\underset{k\rightarrow\infty}{\text{lim}}d\left(x^{k},\text{Fix}T\right)=0.\label{eq:BR-def}
\end{equation}
\end{defn}

In \cite[Proposition 4.1(i)]{CegielskiZalas} it is shown that for
a BR operator $T$, $T-Id$ is demi-closed at $0$, where $Id$ is
the identity operator. Here and throughout, we denote the strong and
weak convergence of a sequence $\left(x^{k}\right)_{k\in\mathbb{N}}$
to a limit point $x\in\mathcal{H}$ by $x^{k}\rightarrow x$ and $x^{k}\rightharpoonup x$,
respectively.
\begin{defn}
\label{def:demi-closed} An operator $T:D\rightarrow\mathcal{H}$
is demi-closed at $0$ if for any sequence $x^{k}\rightharpoonup y\in D$
with $T\left(x^{k}\right)\rightarrow0$ we have $Ty=0$.
\end{defn}

\begin{prop}
\label{prop: BR is DC} \cite[Proposition 4.1(i)]{CegielskiZalas}
Let $T:\mathcal{H}\rightarrow\mathcal{H}$ with $\text{Fix}T\neq\emptyset$
be QNE. If $T$ is BR, then $T-Id$ is demi-closed at $0$.
\end{prop}

The projection operator is a type of a cutter operator. This class
is known in the literature by several alternative names, such as $\mathcal{T}$-class
\cite{BauschkeCombettes_paper}, firmly quasi-nonexpansive \cite[Section B]{quasiShrinking},
directed operators \cite{SegalCensor} and separating operators \cite{CegielskiCutters}.
The name cutter was proposed by Cegielski and Censor in \cite{CeigielskiCensor}.
\begin{defn}
\label{def: cutter}An operator $T:\mathcal{H}\rightarrow\mathcal{H}$
is a cutter if for all $x\in\mathcal{H}$ and $z\in\text{Fix}T\neq\emptyset$,
we have 
\begin{equation}
\left\langle x-T\left(x\right),z-T\left(x\right)\right\rangle \leq0.
\end{equation}

By a similar consideration to the one made previously, it is easy
to verify that Equation (\ref{eq: characterization and Lemma 1.2.5})
can also be derived for a cutter. Nevertheless, it is important to
note that this property also follows from the known fact that a cutter
is a 1-strongly quasi-nonexpansive operator; see, e.g., \cite[Definition 2.3.(iii)]{ReichZalasBarsahd}
and \cite[Page 304]{ReichZalas}.
\end{defn}

We proceed by recalling Opial’s Lemma, a fundamental tool for convergence
analysis in real Hilbert spaces. Since this fact is not fully documented
in the original paper \cite{Opial}, we provide a version adapted
to our needs, based on \cite[Fact 1.1]{ArakcheevBauschke}. 
\begin{lem}
\label{lem: Opial's Lemma} (Opial's Lemma) Let $D$ be a nonempty
closed and convex subset of $\mathcal{H}$ and let $\left(x^{k}\right)_{k\in\mathbb{N}}$
be a sequence in $\mathcal{H}$ such that for all $p\in D$, 
\begin{equation}
\underset{k\rightarrow\infty}{\text{lim}}\left\Vert x^{k}-p\right\Vert 
\end{equation}
exists. If all weak cluster points of $\left(x^{k}\right)_{k\in\mathbb{N}}$
lie in $D$, then $\left(x^{k}\right)_{k\in\mathbb{N}}$ converges
weakly to some $x\in D$.
\end{lem}

Next, the property of bounded regularity for sets is introduced.
\begin{defn}
\label{def:BR sets}Let $\left(D_{i}\right)_{i\in I}$ be a family
of closed and convex sets with a nonempty intersection $D$. We say
that this family is boundedly regular if for any bounded sequence
$\left(x^{k}\right)_{k\in\mathbb{N}}\subseteq\mathcal{H}$ the following
implication holds: 
\begin{equation}
\underset{k\rightarrow\infty}{\text{lim}}\text{\ensuremath{\underset{i\in I}{\text{max}}} }d\left(x^{k},D_{i}\right)=0\:\Longrightarrow\:\underset{k\rightarrow\infty}{\text{lim}}d\left(x^{k},D\right)=0.
\end{equation}

It is worth noting that the fixed point set of a cutter is closed
and convex; see, for example, \cite[Proposition 2.6(ii)]{BauschkeCombettes_paper};
thus, bounded regularity can be assumed on it.
\end{defn}

\section{The SA Algorithmic Scheme}\label{sec:The-string-averaging-algorithmic}

In this section, the Skeletal SA algorithmic scheme is introduced
alongside its core concepts. To demonstrate the structural flexibility
of this framework, several conventional schemes, including our proposed
ROC-SA method, are shown to align with the Skeletal SA structure.
\begin{defn}
\label{def:IndexVec-StingrOp-FitOmega-FiniteCase-1} Let $\left(T_{i}\right)_{i\in I}$
be a finite family of operators on $\mathcal{H}$. Let the ``string''
$I_{t}$ be a finite nonempty ordered subset of $I$, for $t=1,2,\ldots,J$,
of the form $I_{t}=\left(i^{t}_{1},i^{t}_{2},\ldots,i^{t}_{\gamma\left(I_{t}\right)}\right)$,
where the length of the string $I_{t}$, denoted by $\gamma\left(I_{t}\right)$,
is the number of elements in it. Define the operator $T\left[t\right]:\mathcal{H}\rightarrow\mathcal{H}$
as the ``$t$-th string response of the operators $\left(T_{i}\right)_{i\in I_{t}}$'',
and call it a ``string operator'', and the operator $R:\mathcal{H^{J}}\rightarrow\mathcal{H}$
as the ``response of the string operators'', termed the ``combining
operator''. Finally, let $\varOmega:=\left\{ I_{t}\right\} ^{J}_{t=1}$
a set of strings.
\end{defn}

\begin{algorithm}[H]
\caption{\textbf{The Skeletal SA algorithmic scheme }}\label{alg:The-Skeletal-String-Averaging}

Initialization: Choose an arbitrary $x^{0}\in\mathcal{H}$.\\
Iterative Step: Given the current iterate $x^{k}$, choose a set of
$J_{k}$ strings, $\varOmega_{k}$, and \\
(i) calculate, for every $t_{k}=1,2,\ldots,J_{k}$,
\begin{equation}
T\left[t_{k}\right]\left(x^{k}\right),
\end{equation}
(ii) and then calculate 
\begin{equation}
x^{k+1}=R\left(T\left[1\right]\left(x^{k}\right),T\left[2\right]\left(x^{k}\right),\ldots,T\left[J_{k}\right]\left(x^{k}\right)\right).
\end{equation}
\end{algorithm}

\subsection{The original SA scheme }

By defining each string response operator $T\left[t\right]$ to be
a composition of operators along the indices of the string $I_{t}$,
while using a fixed set of strings across all iterations, the original
SA method \cite{Pink-book} is obtained.

\begin{algorithm}[H]
\caption{\textbf{The Original SA scheme }}\label{alg: original SA scheme}
Initialization: Choose an arbitrary $x^{0}\in\mathcal{H}$ and a fixed
set $\varOmega$ of $J$ strings.\\
Iterative Step: Given the current iterate $x^{k}$, \\
(i) Calculate, for every $t=1,2,\ldots,J$, 
\begin{equation}
T\left[t\right]\left(x^{k}\right):=T_{i^{t}_{\gamma\left(I_{t}\right)}}\cdots T_{i^{t}_{2}}T_{i^{t}_{1}}\left(x^{k}\right),
\end{equation}

(ii) and then calculate 
\begin{equation}
x^{k+1}=R\left(T\left[1\right]\left(x^{k}\right),T\left[2\right]\left(x^{k}\right),\ldots,T\left[J\right]\left(x^{k}\right)\right).
\end{equation}
\end{algorithm}

\begin{defn}
\label{def:IndexVec-StingrOp-FitOmega-FiniteCase} Given a positive
weight vector $\omega:=\left(\omega_{t}\right)^{J}_{t=1}$, i.e.,
$\omega_{t}>0$ and $\sum^{J}_{t=1}\omega_{t}=1$, set the combining
operator $R:=\sum^{J}_{t=1}\omega_{t}T\left[t\right]$ (which is commonly
used in the context of SA methods. see, e.g., \cite{BRZ,Zaslavski,CombettesYamada,CQmethods}).
Denote by $\mathcal{M}$ the collection of all pairs\textbf{ $\left(\varOmega,\omega\right)$.}
A finite set of strings $\varOmega$ is called \textit{fit} if for
each $i\in I$, there exists a string $I_{t}=\left(i^{t}_{1},i^{t}_{2},\ldots,i^{t}_{\gamma\left(I_{t}\right)}\right)\in\varOmega$
such that $i^{t}_{\ell}=i$ for some $\ell\in\left\{ 1,2,\ldots,\gamma\left(I_{t}\right)\right\} $.
\end{defn}

Under the setup described in Definition \ref{def:IndexVec-StingrOp-FitOmega-FiniteCase}
and by choosing the pairs $\left(\varOmega,\omega\right)\in\mathcal{M}$
from a pre-defined set $\mathcal{M}^{*}\subset\mathcal{M}$, in which
the length of every $I_{t}\in\varOmega$ is bounded and $\omega_{t}$
is bounded away from zero, the DSA procedure emerges; see, for example,
\cite{CensorZaslavski1,CensorZaslavski2}. Specifically, suppose that
$I=\left\{ 1,2,\ldots,m\right\} $, fix a number $\triangle\in\left(0,\frac{1}{m}\right)$
and an integer $\bar{q}\geq m$. Then 
\begin{equation}
\mathcal{M}_{*}\equiv\mathcal{M}_{*}\left(\triangle,\bar{q}\right):=\,\left\{ \left(\varOmega,\omega\right)\in\mathcal{M}\,|\,\gamma\left(I_{t}\right)\leq\bar{q}\:\text{and}\,\triangle\leq\omega_{t},\,\text{for}\,\text{all}\,I_{t}\in\varOmega\right\} .\label{eq:M*-2-1}
\end{equation}

To account for potential variations at each iteration, let $\left(\varOmega_{k},\omega_{k}\right)\in\mathcal{\mathcal{M}}$
denote the $k$-th ordered pair used in the $k^{\text{th}}$ iteration,
$I_{t_{k}}$ a string in $\varOmega_{k}$ and $\omega_{t_{k}}$ its
corresponding weight, for every $t_{k}=1,2,\ldots,J_{k}$. 

{\bfseries{}
\begin{algorithm}[H]
\caption{\textbf{The DSA algorithmic scheme}}\label{alg:The-Dynamic-String-Averaging}

Initialization: Choose an arbitrary $x^{0}\in\mathcal{H}$.\\
Iterative Step: Given the current iterate $x^{k}$, choose $\left(\varOmega_{k},\omega_{k}\right)\in\mathcal{\mathcal{M}_{*}}$
and \\
(i) Calculate, for every $t=1,2,\ldots,J$, 
\begin{equation}
T\left[t_{k}\right]\left(x^{k}\right)=T_{i^{t_{k}}_{\gamma\left(I_{t_{k}}\right)}}\cdots T_{i^{t_{k}}_{2}}T_{i^{t_{k}}_{1}}\left(x^{k}\right),
\end{equation}
(ii) and then calculate 
\begin{equation}
x^{k+1}=\sum^{J_{k}}_{t_{k}=1}\omega_{t_{k}}T\left[t_{k}\right]\left(x^{k}\right).
\end{equation}
\end{algorithm}
}

\subsection{The component-averaged row projections (CARP) method}

Introduced by Gordon and Gordon \cite{GordonANDGordon}, this approach
provides another structure that agrees with Algorithm \ref{alg:The-Skeletal-String-Averaging}.
It belongs to a class of schemes designed for linear systems, where
the rows of the underlying matrix, representing the equations, are
partitioned into arbitrary groups, called blocks. Within each block,
the current iterate is projected sequentially from one hyperplane
to the next with the aid of Kaczmarz steps with relaxations. Following
these updates, each component of the next iterate is determined by
averaging the corresponding components of the resulting block vectors
among only the active blocks that modified it.

The following definitions are presented to provide a precise overview
of this procedure.
\begin{defn}
\textbf{\label{def:The-KACZ-Sweep} The KACZ Sweep Operator}. Let
$B_{t}$, $1\leq t\leq J$, be a block of equations from the linear
system $Ax=b$, where each equation defines a hyperplane, and let
$\left\{ i_{1},i_{2},\ldots,i_{r}\right\} $ be the set of indices
of $B_{t}$'s equations. For any $x\in\mathbb{R}^{n}$, the KACZ sweep
operator $\text{KSWP}\left(B_{t},x\right):\mathcal{\mathbb{R}}^{n}\rightarrow\mathcal{\mathbb{R}}^{n}$
is defined sequentially by setting $x^{0}=x$, and for each $1\leq j<r$,
\begin{equation}
x^{j+1}=x^{j}+\lambda_{j}\frac{b_{i_{j}}-\left\langle a_{i_{j}},x^{j}\right\rangle }{\left\Vert a_{i_{j}}\right\Vert ^{2}}a_{i_{j}},
\end{equation}
where $\lambda_{j}$ is a relaxation parameter. The final output of
the block sweep is given by $\text{KSWP}\left(B_{t},x\right):=x^{r}$.
\end{defn}

\begin{defn}
\textbf{\label{def:The-Component-Averaging-Operator}The Component-Averaging
Operator.} Let $B := \{ B_1, \allowbreak B_2, \allowbreak \dots, \allowbreak B_J \}$
a collection of blocks. The component-averaging operator $\text{CA}_{B}:\left(\mathcal{\mathbb{R}}^{n}\right)^{J}\rightarrow\mathcal{\mathbb{R}}^{n}$
maps $J$ vectors $x_{1},x_{2},\ldots,x_{J}$ in $\mathcal{\mathbb{R}}^{n}$
to a single point whose $\ell^{\text{th}}$ component is defined by
\begin{equation}
\text{CA}_{B}\left(x_{1},x_{2},\ldots,x_{J}\right)_{\ell}:=\frac{1}{s_{\ell}}\sum_{t\in I_{\ell}}x^{t}_{\ell},
\end{equation}
\end{defn}

where $x^{t}_{\ell}$ is the $\ell^{\text{th }}$component of the
vector $x^{t}$, $I_{\ell}$ is the set of block indices containing
equations that modify the $\ell^{\text{th }}$ component and $s_{\ell}:=\left|I_{\ell}\right|$
is the total number of such blocks. 

The following algorithm is the CARP scheme, written in a similar manner
as in \cite[Algorithm 2]{GordonANDGordon}.

{\bfseries{}
\begin{algorithm}[H]
\caption{\textbf{The CARP algorithmic scheme}}\label{alg:CARP-gordon}

Initialization: $x^{0}:=\left(x^{0}_{1},x^{0}_{2},\ldots,x^{0}_{n}\right)\in\mathcal{\mathbb{R}}^{n}$.\\
Iterative Step: Given the current iterate $x^{k}$, do: \\
(i) For every $1\leq t\leq J$, in parallel, execute some finite number
of KACZ sweeps on all the equations of $B_{t}$. For every $1\leq t\leq J$,
denote the resulting point by $\bar{x}^{t}$. Formally, for some $p_{k}\geq1$
(representing the number of internal sweeps), $\bar{x}^{t}=\text{KSWP}^{p_{k}}\left(B_{t},x^{k}\right)$.\\
(ii) $x^{k+1}=\text{CA}_{B}\left(\bar{x}^{1},\bar{x}^{2},\ldots,\bar{x}^{J}\right)$.
\end{algorithm}
}

To formulate this approach under Definition \ref{def:IndexVec-StingrOp-FitOmega-FiniteCase-1},
fix $p_{k}=1$ across all iterations, set $\varOmega=\left\{ B_{1},B_{2},\ldots,B_{J}\right\} $,
let each string response operator $T\left[t\right]\left(x^{k}\right):=\text{KSWP}^{p_{k}}\left(B_{t},x^{k}\right)$
for every $1\leq t\leq J$, and set $R:=\text{CA}_{B}\left(T\left[1\right]\left(x^{k}\right),T\left[2\right]\left(x^{k}\right),\right.$
\\
$\left.,\ldots,T\left[J\right]\left(x^{k}\right)\right)$.

\subsection{The ROC DSA method }

We adapt the framework of Algorithm \ref{alg:The-Dynamic-String-Averaging}
to form our ROC DSA method. Loosely speaking, given an ordered pair
$\left(\varOmega_{k},\omega_{k}\right)\in\mathcal{M}$, we identify
the operators whose indices appear in it, calculate their distances
from the current iterate, and define the ``farthest'' operator as\textcolor{red}{{}
}the string's end-point. These newly defined end-points are then combined
with a convex combination. By doing so, we create a new class of DSA
methods, and in particular SA ones, where each string response operator
$T\left[t_{k}\right]$ in the Skeletal SA structure corresponds to
selecting the remotest operator along a given string. An important
observation here is that the choice of $\left(\varOmega_{k},\omega_{k}\right)$
is indeed from $\mathcal{M}$ instead of $\mathcal{M}_{*}$, as will
be further detailed in Section \ref{sec:Main-Results}.
\begin{defn}
\label{def:RSC SA} Let $\left(T_{i}\right)_{i\in I}$ be a finite
family of operators, let $k\in\mathbb{N}$ and let $\left(\varOmega_{k},\omega_{k}\right)\in\mathcal{M}$
with a fit finite set $\varOmega_{k}$. For every $t_{k}=1,2,\ldots,J_{k}$,
we denote by
\begin{equation}
i^{t_{k}}_{s}:=\text{min \ensuremath{\left\{  \text{\text{Argmax \ensuremath{\left\{  \left\Vert x^{k}-T_{i^{t_{k}}_{\ell}}\left(x^{k}\right)\right\Vert \mid\ell=1,2,\ldots,\gamma\left(I_{t_{k}}\right)\right\} } }}\right\} } }\label{eq:RSC-control index}
\end{equation}
the ``remotest operator's index of the $t_{k}$-th string'', and
by $\hat{T}_{t_{k}}:=T_{i^{t_{k}}_{s}}$ the remotest operator associated
with it, where the index $s$ itself is $s\equiv s\left(t_{k}\right)$. 
\end{defn}

In the presence of multiple maximizers, the minimum index is chosen
to guarantee uniqueness. For the reader's convenience, we provide
a simple example to clarify this control mechanism.
\begin{example}
\label{exa:ROC-SA} Let $T_{1},T_{2},T_{3}$ be a family of operators
from $\mathbb{R}^{n}$ to $\mathbb{R}$. Given the $k^{\text{th}}$
iterate, assume that $x^{k}=0$ and that $T_{1}\left(0\right)=1$,
$T_{2}\left(0\right)=5$, $T_{3}\left(0\right)=2$. Let $\varOmega_{k}$
consist of a single string $I_{1}:=\left(1,2,3\right)$. Then  $\min \{ \text{Argmax} \{ \| x^k - T_{i_\ell^1} (x^k) \| \mid \ell = 1,2,3 \} \} = \min \{ \text{Argmax} \{ \| T_{i_\ell^1} (0) \| \mid \ell = 1,2,3 \} \} = 2.$
Hence, $i^{t_{k}}_{s}=i^{1}_{2}$ and $\hat{T}_{1}:=T_{i^{1}_{2}}$.\\

We are now ready to introduce our remotest operator control DSA method.
\\
\begin{algorithm}[H]
\caption{\textbf{The ROC DSA method }}\label{alg:The-remotest-operator scheme}

Initialization: Let $\left(T\right)_{i\in I}$ be a finite family
of operators and choose an arbitrary $x^{0}\in\mathcal{H}.$\\
Iterative Step: Given the current iterate $x^{k}$, pick an ordered
pair $\left(\varOmega_{k},\omega_{k}\right)\in\mathcal{M},$ with
a fit finite set $\varOmega_{k},$ and calculate\\
\begin{equation}
x^{k+1}=\sum^{J_{k}}_{t_{k}=1}\omega_{t_{k}}\hat{T}_{t_{k}}\left(x^{k}\right).
\end{equation}
\end{algorithm}
\end{example}

\section{Main Results}\label{sec:Main-Results}

In this section we show that Algorithm \ref{alg:The-remotest-operator scheme}
solves the CFPP for a finite family of BR cutters, $\left(T_{i}\right)_{i\in I}$.
We begin by proving the following lemmata.
\begin{lem}
\label{lem:Fejer-monotone}Any sequence $\left(x^{k}\right)_{k\in\mathbb{N}}$
generated by Algorithm \ref{alg:The-remotest-operator scheme} is
Fej\'{e}r monotone. Namely, for every $z\in F$, 
\begin{equation}
\left\Vert x^{k+1}-z\right\Vert \leq\left\Vert x^{k}-z\right\Vert .
\end{equation}
\end{lem}

\begin{proof}
Let $z\in F$. By Definitions \ref{def:IndexVec-StingrOp-FitOmega-FiniteCase}
and \ref{def:RSC SA}, and the convexity of the function $x\rightarrow\left\Vert x\right\Vert ^{2}$;
see, e.g., the proof in \cite[Lemma 9]{CensorZaslavski1}, it follows
that 
\begin{align}
\left\Vert x^{k+1}-z\right\Vert ^{2}= & \left\Vert \sum^{J_{k}}_{t_{k}=1}\omega_{t_{k}}\hat{T}_{t_{k}}\left(x^{k}\right)-z\right\Vert ^{2}\\
= & \left\Vert \sum^{J_{k}}_{t_{k}=1}\omega_{t_{k}}T_{i^{t_{k}}_{s}}\left(x^{k}\right)-z\right\Vert ^{2}\\
= & \left\Vert \sum^{J_{k}}_{t_{k}=1}\omega_{t_{k}}T_{i^{t_{k}}_{s}}\left(x^{k}\right)-\sum^{J_{k}}_{t_{k}=1}\omega_{t_{k}}z\right\Vert ^{2}\\
\leq & \sum^{J_{k}}_{t_{k}=1}\omega_{t_{k}}\left\Vert T_{i^{t_{k}}_{s}}\left(x^{k}\right)-z\right\Vert ^{2}.
\end{align}
So, by using Inequality (\ref{eq: characterization and Lemma 1.2.5})
for each operator in the family $\left(T_{i^{t_{k}}_{s}}\right)^{J_{k}}_{t_{k}=1}$
and by considering Definition \ref{def:IndexVec-StingrOp-FitOmega-FiniteCase},
\begin{align}
\sum^{J_{k}}_{t_{k}=1}\omega_{t_{k}}\left\Vert T_{i^{t_{k}}_{s}}\left(x^{k}\right)-z\right\Vert ^{2}\leq & \sum^{J_{k}}_{t_{k}=1}\omega_{t_{k}}\left(\left\Vert x^{k}-z\right\Vert ^{2}-\left\Vert T_{i^{t_{k}}_{s}}\left(x^{k}\right)-x^{k}\right\Vert ^{2}\right)\label{eq: for lemma 6-1}\\
= & \left\Vert x^{k}-z\right\Vert ^{2}-\sum^{J_{k}}_{t_{k}=1}\omega_{t_{k}}\left\Vert T_{i^{t_{k}}_{s}}\left(x^{k}\right)-x^{k}\right\Vert ^{2}.\label{eq:for lemma 6-2}
\end{align}
Hence, because of the non-negativity of the norm and weights, we deduce
that 
\begin{equation}
\left\Vert x^{k+1}-z\right\Vert ^{2}\leq\left\Vert x^{k}-z\right\Vert ^{2},
\end{equation}
and the required result is obtained.
\end{proof}

The Fej\'{e}r monotonic behavior of our approach in Lemma \ref{lem:Fejer-monotone}
motivates the following tool.
\begin{lem}
\label{lem: x^k+1-x^k} Assume that for every $k\in\mathbb{N}$ and
for every $t_{k}=1,2,\ldots,J_{k}$, $\omega_{t_{k}}\in\left[\triangle,1\right]$
for some $\triangle>0$. Then any sequence $\left(x^{k}\right)_{k\in\mathbb{N}}$
generated by Algorithm \ref{alg:The-remotest-operator scheme} satisfies
\begin{equation}
\left\Vert x^{k+1}-x^{k}\right\Vert \rightarrow0.
\end{equation}
\end{lem}

\begin{proof}
Equation (\ref{eq:for lemma 6-2}) yields that for every $z\in F$,
\begin{equation}
\left\Vert x^{k+1}-z\right\Vert ^{2}\leq\left\Vert x^{k}-z\right\Vert ^{2}-\sum^{J_{k}}_{t_{k}=1}\omega_{t_{k}}\left\Vert T_{i^{t_{k}}_{s}}\left(x^{k}\right)-x^{k}\right\Vert ^{2}.
\end{equation}
Rearranging the terms, we arrive at
\begin{equation}
\sum^{J_{k}}_{t_{k}=1}\omega_{t_{k}}\left\Vert T_{i^{t_{k}}_{s}}\left(x^{k}\right)-x^{k}\right\Vert ^{2}\leq\left\Vert x^{k}-z\right\Vert ^{2}-\left\Vert x^{k+1}-z\right\Vert ^{2}.\label{eq: x^k-x^k converges to zero - 1}
\end{equation}
Consider the series $\sum_{k\in\mathbb{N}}\left(\sum^{J_{k}}_{t_{k}=1}\omega_{t_{k}}\left\Vert T_{i^{t_{k}}_{s}}\left(x^{k}\right)-x^{k}\right\Vert ^{2}\right)$
and its sequence of partial sums, $\left(\sum^{k}_{r=0}\left(\sum^{J_{r}}_{t_{r}=1}\omega_{t_{r}}\left\Vert T_{i^{t_{r}}_{s}}\left(x^{r}\right)-x^{r}\right\Vert ^{2}\right)\right)_{k\in\mathbb{N}}$.
From (\ref{eq: x^k-x^k converges to zero - 1}), we obtain for any
$z\in F$, recursively, for every $k\in\mathbb{N}$ that
\begin{equation}
\sum^{k}_{r=0}\left(\sum^{J_{r}}_{t_{r}=1}\omega_{t_{r}}\left\Vert T_{i^{t_{r}}_{s}}\left(x^{r}\right)-x^{r}\right\Vert ^{2}\right)\leq\left\Vert x^{0}-z\right\Vert ^{2}-\left\Vert x^{k+1}-z\right\Vert ^{2}\leq\left\Vert x^{0}-z\right\Vert ^{2},
\end{equation}
thus, $\left(\sum^{k}_{r=0}\left(\sum^{J_{r}}_{t_{r}=1}\omega_{t_{r}}\left\Vert T_{i^{t_{r}}_{s}}\left(x^{r}\right)-x^{r}\right\Vert ^{2}\right)\right)_{k\in\mathbb{N}}$
is bounded. Furthermore, it is also an increasing sequence since every
added element is positive, and hence, converges. As a convergent series
and due to our choice of weights, 
\begin{equation}
\sum^{J_{k}}_{t_{k}=1}\omega_{t_{k}}\left\Vert T_{i^{t_{k}}_{s}}\left(x^{k}\right)-x^{k}\right\Vert ^{2}\rightarrow0,\label{eq: x^k-x^k convergence - 2}
\end{equation}
which implies that 
\begin{equation}
\left\Vert T_{i^{t_{k}}_{s}}\left(x^{k}\right)-x^{k}\right\Vert \rightarrow0.\label{eq:x^k-x^k convergence - 3}
\end{equation}
 From Definitions \ref{def:IndexVec-StingrOp-FitOmega-FiniteCase}
and \ref{def:RSC SA} it follows that 
\begin{align}
\left\Vert x^{k+1}-x^{k}\right\Vert  & =\left\Vert \sum^{J_{k}}_{t_{k}=1}\omega_{t_{k}}\left(T_{i^{t_{k}}_{s}}\left(x^{k}\right)-x^{k}\right)\right\Vert \\
 & \leq\sum^{J_{k}}_{t_{k}=1}\omega_{t_{k}}\left\Vert T_{i^{t_{k}}_{s}}\left(x^{k}\right)-x^{k}\right\Vert ,
\end{align}
and so, because of (\ref{eq:x^k-x^k convergence - 3}), $\left\Vert x^{k+1}-x^{k}\right\Vert \rightarrow0$.
\end{proof}

Following Bauschke and Combettes' work \cite{BauschkeCombettes_paper},
the authors of \cite{ReichZalasBarsahd} developed a displacement
method \cite[Lemma 4.2]{ReichZalasBarsahd}, based on the remotest
set control principles. Their approach makes use of the so-called
$p$-intermittent sequence, see, for example, \cite[Definition 5.8.8]{Cegielski},
which guarantees that every constraint index appears infinitely many
times in the iterative procedure.
\begin{defn}
A sequence $\left(I_{k}\right)_{k\in\mathbb{N}}\subseteq I$ is a
$p$-intermittent sequence, where $p\in\mathbb{N}$, if $I=I_{1}\uplus I_{2}\uplus\ldots\uplus I_{k+p-1}$
for any $k\in\mathbb{N}$.
\end{defn}

The $p$-intermittent sequence plays a key role in our main result.
Instead of being directly assumed, this property arises naturally
from our use of a repetitive sequence, also known as an an admissible
sequence, see, for example, \cite[Definition 5.6.11]{Cegielski},
that encompasses all constraint indices. Specifically, for every $i\in I$,
we assume that there exists an $N_{i}\in\mathbb{N}$ such that for
every $k\in\mathbb{N}$, $i\in\bigcup^{k+N_{i}-1}_{m=k}I_{m}$. By
letting $p:=\text{max}_{i\in I}N_{i}$, the sequence $\left(I_{k}\right)_{k\in\mathbb{N}}$
is inherently a $p$-intermittent one. While this structural condition
is mild, allowing the repetition interval $M_{i}$ to depend entirely
on the individual index $i$ significantly relaxes the standard, strict
requirements of $p$-intermittency. Under these relaxed considerations,
we establish the convergence of Algorithm \ref{alg:The-remotest-operator scheme}.
\begin{thm}
\label{thm:main result}Let $\left(T_{i}\right)_{i\in I}$ be a finite
family of BR cutters such that $F:=\cap_{i\in I}\text{Fix}T_{i}\neq\emptyset$,
let $\left(\varOmega_{k},\omega_{k}\right)\in\mathcal{\mathcal{M}}$
be with a fit finite set $\varOmega_{k}$ for every $k\in\mathbb{N}$
and weights as in Lemma \ref{lem: x^k+1-x^k}, and assume that $I$
is repetitive. Suppose also that the family of sets $\left(\text{Fix}T_{i}\right)_{i\in I}$
is boundedly regular. Then any sequence $\left(x^{k}\right)_{k\in\mathbb{N}}$
generated by Algorithm \ref{alg:The-remotest-operator scheme} converges
strongly to $x^{*}\in F$.
\end{thm}

\begin{proof}
For every $i\in I$ there exists a sequence $\left(r_{i}\left(k\right)\right)_{k\in\mathbb{N}}$
for which, for every $k\in\mathbb{N}$, $r_{i}\left(k\right)\in\left\{ 1,2,\ldots,p\right\} $
and the iterate $x^{k+r_{i}\left(k\right)}$ ``contains'' the $i^{\text{th}}$
index, namely, $i\in I_{k+r_{i}\left(k\right)}$. Hence, by repeatedly
applying the Fej\'{e}r monotonocity of $\left(x^{k}\right)_{k\in\mathbb{N}}$
for $p$-consecutive steps, and by the fitness of $\varOmega$, 
\begin{equation}
\left\Vert x^{k+p}-z\right\Vert ^{2}\leq\left\Vert x^{k}-z\right\Vert ^{2}-\sum_{i\in I}\omega_{t_{k+r_{i}\left(k\right)}}\left\Vert T_{i}\left(x^{k+r_{i}\left(k\right)}\right)-x^{k+r_{i}\left(k\right)}\right\Vert ^{2},\label{eq:main results-1}
\end{equation}
for every $z\in F$. It should be emphasized that in case there exists
an index $j\in I$ for which $r_{i}\left(\bar{k}\right)=j$ and $r_{i}\left(\tilde{k}\right)=j$
with $\bar{k}>\tilde{k}$, then (\ref{eq:main results-1}) follows
because of the norm non-negativity. In this case, we choose an arbitrary
single position where the $j^{\text{th}}$ index appears.

From similar considerations to those made in the proof of Lemma \ref{lem: x^k+1-x^k},
it follows that for every $i\in I$, 
\begin{equation}
T_{i}\left(x^{k+r_{i}\left(k\right)}\right)-x^{k+r_{i}\left(k\right)}\rightarrow0.\label{equ: main results-2}
\end{equation}
We proceed by showing that $\left\Vert x^{k+r}-x^{k}\right\Vert \rightarrow0$
for every $1\leq r\leq p$. By the triangle inequality
\begin{align}
\left\Vert x^{k+r}-x^{k}\right\Vert  & =\left\Vert x^{k+r}-x^{k}+x^{k+1}-x^{k+1}+x^{k+2}-\right.\\
 & \left.-x^{k+2}+\ldots+x^{k+r-1}-x^{k+r-1}\right\Vert \\
 & \leq\sum^{p-1}_{l=0}\left\Vert x^{k+l+1}-x^{k+l}\right\Vert ,
\end{align}
and hence,
\begin{equation}
0\leq\underset{k\rightarrow\infty}{\text{lim}}\left\Vert x^{k+r}-x^{k}\right\Vert \leq\underset{k\rightarrow\infty}{\text{lim}}\sum^{p-1}_{l=0}\left\Vert x^{k+l+1}-x^{k+l}\right\Vert .
\end{equation}
Combining this inequality with Lemma \ref{lem: x^k+1-x^k}, we see
that 
\begin{equation}
\left\Vert x^{k+r}-x^{k}\right\Vert \rightarrow0,\label{eq: main result-3}
\end{equation}
 for every $1\leq r\leq p$.

Now, from the boundeness of $\left(x^{k}\right)_{k\in\mathbb{N}}$
there exists a weakly convergent subsequence $\left(x^{k_{n}}\right)_{n\in\mathbb{N}}$
such that $x^{k_{n}}\rightharpoonup x\in\mathcal{H}$. Let $i\in I$.
From previous discussion, there exists a subsequence of $\left(r_{i}\left(k\right)\right)_{k\in\mathbb{N}}$,
say $\left(r_{i}\left(k_{n}\right)\right)_{n\in\mathbb{N}}$, for
which $i\in I_{k_{n}+r_{i}\left(k_{n}\right)}$ for every $n\in\mathbb{N}$.
Hence, due to (\ref{eq: main result-3}), it follows that 
\begin{equation}
x^{k_{n}}-x^{k_{n}+r_{i}\left(k_{n}\right)}\rightarrow0,\label{eq: main result-4}
\end{equation}
and so, $x^{k_{n}+r_{i}\left(k_{n}\right)}\rightharpoonup x$, as
well. This, combined with Proposition \ref{prop: BR is DC} and $T_{i}\left(x^{k_{n}+r_{i}\left(k_{n}\right)}\right)-x^{k_{n}+r_{i}\left(k_{n}\right)}\rightarrow0$,
which is derived from (\ref{equ: main results-2}), results in $x\in\text{Fix}T_{i}$.
This holds for each $i\in I$, and thus, $x\in F$. Therefore, Opial's
Lemma conditions are satisfied and as a result, $x^{k}\rightharpoonup x^{*}\in F$.
To establish strong convergence note that Equations (\ref{eq:BR-def})
and (\ref{equ: main results-2}) yield that for every $i\in I$, 
\begin{equation}
\underset{n\rightarrow\infty}{\text{lim}}d\left(x^{k+r_{i}\left(k\right)},\text{Fix}T_{i}\right)=0.\label{eq: 3.25}
\end{equation}
By (\ref{eq: inf-distance}) and the nonexpansivity of the distance
function; see, e.g., \cite[Example 1.48]{BauschkeCombettes}, it follows
that for every $i\in I$, 
\begin{equation}
\left\Vert d\left(x^{k+r_{i}\left(k\right)},\text{Fix}T_{i}\right)-d\left(x^{k},\text{Fix}T_{i}\right)\right\Vert \leq\left\Vert x^{k+r_{i}\left(k\right)}-x^{k}\right\Vert ,
\end{equation}
and consequently, because of (\ref{eq: main result-3}) and the sandwich
rule, 
\begin{equation}
d\left(x^{k+r_{i}\left(k\right)},\text{Fix}T_{i}\right)-d\left(x^{k},\text{Fix}T_{i}\right)\rightarrow0,
\end{equation}
which implies that $d\left(x^{k},\text{Fix}T_{i}\right)\rightarrow0$
with the aid of (\ref{eq: 3.25}). Hence, by the bounded regularity
property of $\left(\text{Fix}T_{i}\right)_{i\in I}$, we are able
to deduce that $d\left(x^{k},F\right)\rightarrow0$. Let $\varepsilon>0$.
Then there exists an integer $N_{\varepsilon}\in\mathbb{N}$ and a
$\bar{z}\in F$ for which, $\text{inf}_{z\in F}\left\Vert x^{N_{\varepsilon}}-z\right\Vert =\left\Vert x^{N_{\varepsilon}}-\bar{z}\right\Vert <\frac{\varepsilon}{2}$.
Let $m,q\geq N_{\varepsilon}$. Then by the triangle inequality and
by Lemma \ref{lem:Fejer-monotone}, it follows that 
\begin{equation}
\left\Vert x^{m}-x^{q}\right\Vert \leq\left\Vert x^{m}-\bar{z}\right\Vert +\left\Vert x^{q}-\bar{z}\right\Vert \leq\left\Vert x^{N_{\varepsilon}}-\bar{z}\right\Vert +\left\Vert x^{N_{\varepsilon}}-\bar{z}\right\Vert <\frac{\varepsilon}{2}+\frac{\varepsilon}{2}<\varepsilon.
\end{equation}
As a result, $\left(x^{k}\right)_{k\in\mathbb{N}}$ is a Cauchy sequence,
thus, converges strongly (because of the completeness of $\mathcal{H}$)
to $x^{*}$. Note that the latter holds since every strongly convergent
sequence converges weakly, and weak limits are unique.
\end{proof}

\begin{rem}
\label{rem:infinite strings} The length of the strings in Definition
\ref{def:IndexVec-StingrOp-FitOmega-FiniteCase-1} is restricted to
a finite cardinality. This ensures that every string operator is well
defined. However, due to the ROC properties, in the ROC DSA method
we are able to extend the lengths of the strings to be infinite as
long as $I$ is of a finite cardinality. In this case an infinite
string must have repetitive indices. Therefore, for a given infinite
string $I_{t}$, finding the remotest operator control index of the
$t$-th string $i^{t}_{s}$ is equivalent to finding it for the same
string without repetitions of indices. This insight gives rise to
an additional property that has yet to be studied in the context of
SA methods: Strings of an infinite length. It is worth highlighting
that the set $\mathcal{M}$ encompasses strings of an infinite lengths;
nevertheless, the SA framework only makes use of ordered pairs $\left(\varOmega,\omega\right)\in\mathcal{M}$,
for which, for every $t=1,2,\ldots,J$, $\gamma\left(I_{t}\right)<\infty$.
In our ROC DSA scheme we do not restrict ourselves to this scenario
only.
\end{rem}

It should be emphasized that in case $\text{dim\ensuremath{\mathcal{H}}}<\infty$,
the use of Opial's Lemma is unnecessary and the bounded regularity
assumption on $\left(\text{Fix}T_{i}\right)_{i\in I}$, which is known
to hold inherently in finite dimensional spaces; see, e.g., \cite[Proposition 4.9]{CegielskiZalas},
is not required for our argument. This is so because of the Bolzano-Weierstrass
Theorem, ensuring that the boundeness of $\left(x^{k}\right)_{k\in\mathbb{N}}$
results in a strongly convergent subsequence, rather than a weakly
convergent one, and by the closedness of $\text{Fix}T_{i}$ for every
$i\in I$, that guarantees that $x\in F$. We show this in our alternative
proof of Theorem \ref{thm:main result}, given below, for the finite
dimensional case.
\begin{thm}
\label{thm:special case of main result-1}Assume that $\text{dim}\mathcal{H}<\infty$
and that all the conditions of Theorem \ref{thm:main result}, excluding
the bounded regularity of $\left(\text{Fix}T_{i}\right)_{i\in I}$,
hold. Then any sequence $\left(x^{k}\right)_{k\in\mathbb{N}}$ generated
by Algorithm \ref{alg:The-remotest-operator scheme} converges strongly
to $x^{*}\in F$.
\end{thm}

\begin{proof}
Suppose that $\left(x^{k_{n}}\right)_{n\in\mathbb{N}}\rightarrow\bar{x}$.
Because of the bounded regularity of each of the $T_{i}$-s and Equation
(\ref{eq: main result-3}), $\underset{n\rightarrow\infty}{\text{lim}}d\left(x^{k_{n}+r_{i}\left(k_{n}\right)},\text{Fix}T_{i}\right)=0$,
for every $i\in I$. Additionally, the continuity of the distance
function with Equation (\ref{eq: main result-4}) yield that 
\begin{equation}
\underset{n\rightarrow\infty}{\text{lim}}d\left(x^{k_{n}+r_{i}\left(k_{n}\right)},\text{Fix}T_{i}\right)=d\left(\underset{n\rightarrow\infty}{\text{lim}}x^{k_{n}+r_{i}\left(k_{n}\right)},\text{Fix}T_{i}\right)=d\left(x,\text{Fix}T_{i}\right)=0.
\end{equation}
Hence, due to the closedness of $\text{Fix}T_{i}$ for every $i\in I$,
it follows that $\bar{x}\in F$. Now, since $\left(x^{k}\right)_{k\in\mathbb{N}}$
is Fej\'{e}r monotone, the sequence $\left(\left\Vert x^{k}-\bar{x}\right\Vert \right)_{k\in\mathbb{N}}$
is monotonically-decreasing and bounded, and has a subsequence, $\left(\left\Vert x^{k_{n}}-\bar{x}\right\Vert \right)_{n\in\mathbb{N}},$
which converges to $0$. Consequently, it converges to $0$, as well.
Thus, $x^{k}\rightarrow\bar{x}$.
\end{proof}

A special case of our DSA procedure occurs when a single, pre-determined
ordered pair is used throughout the iterative process in Algorithm
\ref{alg:The-remotest-operator scheme}. This type of scheme is generally
known as \textquotedbl static'', see, e.g, \cite{CensorNisenbaum}.
If one seeks to derive it using the proof methodology of Theorem \ref{thm:main result},
then the repetitiveness assumption on $I$ still holds. Nevertheless,
due to the properties of ROC, this condition is no longer required,
as demonstrated in Theorem \ref{thm: main static result} below. Intuitively,
the idea is that if there exists an operator which is never chosen
as the ``farthest'' one during the iterative procedure, then as
the overall sequence of distances shrinks, the distance to that operator
must also decrease.\\
\begin{algorithm}[H]

\caption{\textbf{The static ROC method}}\label{alg:The-static-ROC}

Initialization: Choose an arbitrary $x^{0}\in\mathcal{H}$ and an
ordered pair $\left(\varOmega,\omega\right)\in\mathcal{M}$ where
$\varOmega$ is a fit finite set.\\
Iterative Step: Given the current iterate $x^{k}$, calculate\\
\begin{equation}
x^{k+1}=\sum^{J}_{t=1}\omega_{t}\hat{T}_{t}\left(x^{k}\right).
\end{equation}
\end{algorithm}

\begin{defn}
For every $x\in\mathcal{H}$ and $t=1,2,\ldots,J$, define the distance
to the remotest operator in the string by 
\begin{equation}
G_{t}\left(x\right):=\text{\ensuremath{\text{max}\left\{ \left\Vert T_{\ell}\left(x\right)-x\right\Vert \mid\ell=1,2,\ldots,\gamma\left(I_{t}\right)\right\} .} }
\end{equation}
\end{defn}

\begin{lem}
\label{lem:sequences converges to zero} Let $\left(\varOmega,\omega\right)\in\mathcal{M}$
be such that for every $t=1,2,\ldots,J$, $\omega_{t}\in\left[\triangle,1\right]$
for some $\triangle>0$. Then for any sequence $\left(x^{k}\right)_{k\in\mathbb{N}}$
generated by Algorithm \ref{alg:The-static-ROC}, we have that for
every $t=1,2,\ldots,J$, the sequence $\left(G_{t}\left(x^{k}\right)\right)_{k\in\mathbb{N}}$
converges to zero.
\end{lem}

\begin{proof}
From similar arguments made in Lemma \ref{lem: x^k+1-x^k}. 
\end{proof}

\begin{thm}
\label{thm: main static result} Let $\left(T_{i}\right)_{i\in I}$
be a finite family of BR cutters such that $F:=\cap_{i\in I}\text{Fix}T_{i}\neq\emptyset$
and let $\left(\varOmega,\omega\right)\in\mathcal{M}$ be with a fit
finite set $\varOmega$, and such that for every $t=1,2,\ldots,J$,
$\omega_{t}\in\left[\triangle,1\right]$ for some $\triangle>0$.
Assume also that the family of sets $\left(\text{Fix}T_{i}\right)_{i\in I}$
is boundedly regular. Then any sequence $\left(x^{k}\right)_{k\in\mathbb{N}}$
generated by Algorithm \ref{alg:The-static-ROC} converges strongly
to $x^{*}\in F$.
\end{thm}

\begin{proof}
Let us pick a $t\in\left\{ 1,2,\ldots,J\right\} $. From Definition
\ref{def:RSC SA}, it follows that for every $1\leq\ell\leq\gamma\left(I_{t}\right)$,

\begin{equation}
0\leq\left\Vert T_{\ell}\left(x^{k}\right)-x^{k}\right\Vert \leq G_{t}\left(x^{k}\right).\label{eq: Theorem7-1}
\end{equation}
Consequently, by Lemma \ref{lem:sequences converges to zero} and
the sandwich rule, 
\begin{equation}
\left\Vert T_{\ell}\left(x^{k}\right)-x^{k}\right\Vert \rightarrow0,\label{eq:sandwich rule}
\end{equation}
for every $1\leq\ell\leq\gamma\left(I_{t}\right)$, and so, 
\begin{equation}
T_{\ell}\left(x^{k}\right)-x^{k}\rightarrow0,\label{eq:Theorem 7-3}
\end{equation}
for every $1\leq\ell\leq\gamma\left(I_{t}\right)$. Since every Fej\'{e}r
monotone sequence is bounded, Lemma \ref{lem:Fejer-monotone} implies
that $\left(x^{k}\right)_{k\in\mathbb{N}}$ is a bounded sequence.
Thus, it has a weakly convergent subsequence, $\left(x^{k_{n}}\right)_{n\in\mathbb{N}}$,
such that 
\begin{equation}
x^{k_{n}}\rightharpoonup\bar{x}\in\mathcal{H}.\label{eq: Theorem-7-5}
\end{equation}
Applying Proposition \ref{prop: BR is DC} on each of the $T_{\ell}$-s
results in the demi-closedness of $T_{\ell}-Id$. Hence, (\ref{eq: Theorem-7-5})
and (\ref{eq:Theorem 7-3}), with $x^{k_{n}}$ in (\ref{eq:Theorem 7-3})
instead of $x^{k}$, results in $x\in\cap^{\gamma\left(I_{t}\right)}_{\ell=1}\text{Fix}T_{\ell}$.
Since our choice of $t$ is arbitrary, this result holds for every
$t=1,2,\ldots,J$, and, therefore, since $\varOmega$ is fit, it follows
that 
\begin{equation}
\bar{x}\in\cap_{i\in I}\text{Fix}T_{i}.\label{eq:Theorem-7-4}
\end{equation}
Fej\'{e}r monotonicity from Lemma \ref{lem:Fejer-monotone} and (\ref{eq:Theorem-7-4})
fulfill Opial's Lemma (Lemma \ref{lem: Opial's Lemma} above) assumptions,
and thus, we deduce that $x^{k}\rightharpoonup x^{*}\in F$.

Equations (\ref{eq:BR-def}) and (\ref{eq:sandwich rule}) imply that
$\underset{k}{\text{lim}\,}d\left(x^{k},\text{Fix}T_{\ell}\right)=0$
for every $1\leq\ell\leq\gamma\left(I_{t}\right)$, for every $t=1,2,\ldots,J$,
and so, by the fitness of $\varOmega$ and the bounded regularity
of $\left(\text{Fix}T_{i}\right)_{i\in I}$, $\underset{k}{\text{lim}\,}d\left(x^{k},F\right)=0$.
Following the final arguments in the proof of Theorem \ref{thm:main result}
results in $x^{k}\rightarrow x^{*}$.

\end{proof}

\begin{rem}
Similarly to the discussion following Remark \ref{rem:infinite strings},
for finite-dimensional spaces, the use of Opial's Lemma becomes unnecessary,
and the bounded regularity assumption on $\left(\text{Fix}T_{i}\right)_{i\in I}$
can be omitted.

\end{rem}

If one takes the family of projections $\left(P_{i}\right)_{i\in I}$
to be in the role of the family $\left(T_{i}\right)_{i\in I}$ in
Algorithm \ref{alg:The-static-ROC}, then the latter turns to be the
remotest set control projection SA method (RSC PSA). Apart from creating
a new SA projection (SAP) algorithmic structure; see, e.g., \cite{CensorTom},
its scheme enables to grasp a better understanding of our general
method by reducing the most violated constraint to the most distant
set. In addition to that, it shows that the well-known greedy projection
algorithm, see, for example, \cite[page 210 and the remark given before Section 1.5]{BeckBook},
as well as the simultaneous projection method; see, e.g., \cite[Section 4.4]{Cegielski}
and the Introduction in \cite{CensorElving}, are special cases of
it, in case $\text{dim\ensuremath{\mathcal{H}}}<\infty$.\textbf{}\\
{\bfseries{}
\begin{algorithm}[H]

\caption{\textbf{RSC-PSA scheme}}\label{alg:RSC-PSA-scheme}

Initialization: Choose an arbitrary $x^{0}\in\mathcal{H}$ and an
ordered pair $\left(\varOmega,\omega\right)\in\mathcal{M}$ where
$\varOmega$ is a fit finite set.\\
Iterative Step: Given the current iterate $x^{k}$, calculate\\
\begin{equation}
x^{k+1}=\sum^{J}_{t=1}\omega_{t}\hat{P}_{t}\left(x^{k}\right).
\end{equation}
\end{algorithm}
}{\bfseries{}
\begin{algorithm}[H]
\caption{\textbf{The greedy projection algorithm}}\label{alg:The-greedy-projection}

Initialization: Choose an arbitrary $x^{0}\in\mathcal{H}$.\\
Iterative Step: Given the current iterate $x^{k}$, calculate\\
\begin{equation}
x^{k+1}=P_{i_{k}}\left(x^{k}\right),
\end{equation}
where the control sequence $\left(i_{k}\right)_{k\in\mathbb{N}}$
is a remotest set one with respect to $I:=\left\{ 1,2,\ldots,m\right\} $.
\end{algorithm}
}
\begin{cor}
Algorithm \ref{alg:The-greedy-projection} is a special case of Algorithm
\ref{alg:RSC-PSA-scheme}.
\end{cor}

\begin{proof}
Suppose that $I=\left\{ 1,2,\ldots,m\right\} $. Take $\varOmega=\left\{ \left(1,2,\ldots,m\right)\right\} $
and $\omega_{1}=1$. Then apply Theorem \ref{thm: main static result}.
\end{proof}

{\bfseries{}
\begin{algorithm}[H]
\caption{\textbf{The simultaneous projection algorithm}}\label{alg:The-simultaneous-projection}

Initialization: Choose an arbitrary $x^{0}\in\mathcal{H}$.\\
Iterative Step: Given the current iterate $x^{k}$, calculate\\
\begin{equation}
x^{k+1}=\sum_{i\in I}\omega_{i}P_{i}\left(x^{k}\right),
\end{equation}
where the $\omega_{i}$-s have similar properties to the weights defined
in Definition \ref{def:IndexVec-StingrOp-FitOmega-FiniteCase}.
\end{algorithm}
}
\begin{cor}
Algorithm \ref{alg:The-simultaneous-projection} is a special case
of Algorithm \ref{alg:RSC-PSA-scheme}.
\end{cor}

\begin{proof}
Suppose that $I=\left\{ 1,2,\ldots,m\right\} $. Take $\varOmega=\left\{ \left(1\right),\left(2\right),\ldots,\left(m\right)\right\} $
and the weights $\omega_{t}=1$, for every $t=1,2,\ldots,J$. Then
apply Theorem \ref{thm: main static result}.
\end{proof}

\section{Concluding Comments}\label{sec:Concluding-Comments}

To the best of our knowledge, our study is the first to investigate
ROC-SA methods. To assist in formulating these methods, a novel structural
framework termed the Skeletal SA framework is presented herein to
allow a sequential step other than compositions. Our greedy-type methods
solve the CFPP as well as its equivalent form, the convex feasibility
problem (CFP); see, e.g., \cite[Section 1.3.4]{Cegielski}, which
is of a great importance too, see, for example, \cite{Censor1984}
and references therein. Furthermore, in the static SA case these methods
converge even if not all of the constraint indices appear in their
iterative procedure. This is in contrast to RSC methods where convergence
cannot be assured unless the latter holds, making this distinction
significant. In addition to that, in both of our static and dynamic
cases, it is possible to use strings of an infinite length, which
have not been considered in the context of the SA scheme. As far as
we are aware, SA schemes utilizing parallel steps other than convex
combinations have not yet been explored. This raises another question
for future research regarding SA methods.

\section*{Statements and Declarations}

\subsection*{Competing Interests}

The authors have no relevant financial or non-financial interests
to disclose.

\subsection*{Use of AI in Research Methodology and Writing}

During the development of the proofs in this paper, the authors utilized
an AI tool to assist in exploring potential proof pathways. This involved
a cyclical process of collaborative brainstorming, where the reasoning
steps co-developed with the tool were critically examined, refined,
and adopted through rigorous conceptual discussion. The final validation
of all mathematical results was conducted entirely by the human authors,
who accept full accountability for their correctness.

\subsection*{Data Availability}

Data sharing is not applicable to this article as no data sets were
generated or analyzed during the current study.

\subsection*{Author Contributions}

All authors jointly developed the theoretical results and approved
the final manuscript.

\section*{Acknowledgments}

The work of Yair Censor was supported by the Cooperation Program between
the German Cancer Research (DKFZ) and Israel's Ministry of Innovation,
Science and Technology (MOST). \bibliography{ROC-SA-CNR-bib}

\end{document}